\documentclass{amsart}
\usepackage{amssymb,amsrefs}
\newtheorem{thm}{Theorem}[section]
\newtheorem{cor}[thm]{Corollary}
\newtheorem{lem}[thm]{Lemma}
\newtheorem{prop}[thm]{Proposition}
\theoremstyle{definition}
\newtheorem{defn}[thm]{Definition}
\theoremstyle{remark}
\newtheorem{rem}[thm]{Remark}
\newtheorem{exa}[thm]{Example}

\newcommand{\K}{\mathbb C}
\newcommand{\Kstar}{\K^\times}
\begin{document}

\title[Gaudin subalgebras and stable rational curves]{Gaudin subalgebras and stable rational curves}%
\author{Leonardo Aguirre}%
\author{Giovanni Felder}%
\author{Alexander P. Veselov}%

\address{Department of mathematics, ETH Zurich, 8092 Zurich, Switzerland}%
\email{giovanni.felder@math.ethz.ch}
\address{Department of mathematics, ETH Zurich, 8092 Zurich, Switzerland}%
\email{leonardo.aguirre@math.ethz.ch}
\address{
Department of Mathematical Sciences,
Loughborough University,
Loughborough
LE11 3TU,
UK and Moscow State University, Moscow 119899, Russia}
\email{A.P.Veselov@lboro.ac.uk}

\thanks{}%

\begin{abstract}
Gaudin subalgebras are abelian Lie subalgebras of maximal dimension
spanned by generators of the Kohno--Drinfeld Lie algebra $\mathfrak
t_n$. We show that Gaudin subalgebras form a variety isomorphic to
the moduli space $\bar M_{0,n+1}$ of stable curves of genus zero
with $n+1$ marked points. In particular, this gives an embedding of
$\bar M_{0,n+1}$ in a Grassmannian of $(n-1)$-planes in an
$n(n-1)/2$-dimensional space. We show that the sheaf of Gaudin
subalgebras over $\bar M_{0,n+1}$ is isomorphic to a sheaf of
twisted first order differential operators. For each representation
of the Kohno--Drinfeld Lie algebra with fixed central character, we
obtain a sheaf of commutative algebras whose spectrum is a
coisotropic subscheme of a twisted version of the logarithmic
cotangent bundle of $\bar M_{0,n+1}$.
\end{abstract}
\maketitle
\section{Introduction}
The Kohno--Drinfeld Lie algebra $\mathfrak t_n$ ($n=2,3,\dots$) over
$\K$, see \cites{Kohno,Drinfeld}, is the quotient of the free Lie
algebra on generators $t_{ij}=t_{ji}$, $i\neq j\in\{1,\dots,n\}$ by
the ideal generated by the relations
\begin{eqnarray*}
{}[t_{ij},t_{kl}]&=&0, \qquad \text{if $i,j,k,l$ are distinct,}\\
{} [t_{ij},t_{ik}+t_{jk}]&=&0, \qquad  \text{if $i,j,k$ are
  distinct.}
\end{eqnarray*}
This Lie algebra appears in \cite{Kohno} as the holonomy Lie algebra
of the complement of the union of the diagonals $z_i=z_j, i<j$ in
$\K^n$. The universal Knizhnik--Zamolodchikov connection
\cite{Drinfeld} takes values in $\mathfrak t_n$.

In this paper we consider the abelian Lie subalgebras of maximal
dimension contained in the linear span $\mathfrak t^1_n$ of the
generators $t_{ij}$. Motivating examples are the algebras considered
by Gaudin \cites{Gaudin1976,Gaudin} in the framework of integrable
spin chains in quantum statistical mechanics and the Jucys--Murphy
subalgebras spanned by
$t_{12},t_{13}+t_{23},t_{14}+t_{24}+t_{34},\dots$, appearing in the
representation theory of the symmetric group (see
\cites{VershikOkounkov1996,VershikOkounkov} and references therein).

Our main result is the classification of Gaudin subalgebras. We show
that they are parametrised by the moduli space $\bar M_{0,n+1}$ of
stable curves of genus zero with $n+1$ marked points (Theorem
\ref{t-1}). The Gaudin subalgebras parametrised by the open subset
$M_{0,n+1}$ are the ones considered originally by Gaudin (with
$t_{ij}$ replaced by their image in certain representations of
$\mathfrak t_n$.) To prove this theorem it is useful to represent
$\bar M_{0,n+1}$ as a subvariety of a product of projective lines
given by explicit equations. We give such a description, proving a
variant of a theorem of Gerritzen, Herrlich and van der Put
\cite{Gerritzenetal}, in the Appendix.

Gaudin subalgebras form a locally free sheaf of Lie algebras on
$\bar M_{0,n+1}$. We describe this sheaf as a sheaf of first order
twisted logarithmic differential operators (Theorem \ref{t-102}.)
For an algebra homomorphism $U\mathfrak t_n\to A$ from the universal
enveloping algebra of $t_n$ to an associative algebra $A$, we then
get a sheaf of commutative subalgebras $\mathcal E_A$ of the
$\mathcal O_X$-algebra $A\otimes \mathcal O_X$ on $X=\bar
M_{0,n+1}$. We show that its relative spectrum is a coisotropic
subscheme of a Poisson variety, a twisted version of the logarithmic
cotangent bundle of $\bar M_{0,n+1}$ (Corollary \ref{c-104}.) For a
large class of representations of $U\mathfrak t_n$ these spectra, or
at least their part over $M_{0,n+1}$, have been recently described
in algebro-geometric terms using the Bethe ansatz, see
\cites{Frenkel,FeiginFrenkelRybnikov, FeiginFrenkelToledanoLaredo,MTV1,MTV2,MTV3} and references therein, and shown to have
surprising connection with several other mathematical subjects. It
will be interesting to relate these descriptions to the geometry of
$\bar M_{0,n+1}$. This and possible generalisations to other root
systems will be the subject of further investigation.

It is interesting to look at our result in the context of the
relation between flag varieties and configuration spaces initiated
by Atiyah \cites{A1, A2, AB}, who was inspired by Berry and Robbins.
Note that the usual flag variety $U(n)/T^n$ can be naturally viewed
as the space of all Cartan subalgebras in the unitary Lie algebra
$\mathfrak u(n).$ Our result and a parallel between Gaudin and
Cartan subalgebras give another link between these two varieties.

The limiting behaviour of the algebras introduced by Gaudin, that in
our approach are parametrised by the open subset $M_{0,n+1}$, has
been studied in various contexts. In \cite{Vinberg} Vinberg studied
the commutative subspaces of degree 2 of the universal enveloping algebra $U\mathfrak g$ of a semisimple Lie algebra $\mathfrak
g$ in relation with Poisson commutative subalgebras of the Poisson algebra $S\mathfrak
g$ of polynomial functions on the dual of $\mathfrak
g$.
In the case of $sl_n$ his result implies a set-theoretic description
of all possible limits of Gaudin subalgebras. Limits of Gaudin
subalgebras of $U(\mathfrak g)^{\otimes n}$ were studied more
recently in \cites{ChervovFalquiRybnikov1,ChervovFalquiRybnikov2}.
In \cite{ChervovFalquiRybnikov1} it is noticed that Jucys-Murphy
elements arise as limits of Gaudin Hamiltonians, see Remark
\ref{r-1}.

It is important to mention that our result works over any field. In
particular, it holds over reals, which is important for
applications. It is known that the set of real points $\bar
M_{0,n+1}(\mathbb R) \subset \bar M_{0,n+1}$ is a smooth real
manifold, which can be glued of $n!/2$ copies of the Stasheff
associahedron (see \cite{Kap}). This gives a very convenient
geometric representation of all limiting cases of  the real Gaudin
subalgebras and related quantum integrable systems. In particular,
Jucys-Murphy subalgebra corresponds to one of the vertices of the
associahedron. The spectrum in this case was studied in detail by
Vershik and Okounkov \cites{VershikOkounkov1996,VershikOkounkov}.
What happens at other vertices labelled by different triangulations
of an $n$-gon is worthy to investigate further. As it was explained
in \cites{FalquiMusso,
ChervovFalquiRybnikov1,ChervovFalquiRybnikov2} the corresponding
integrable systems have a nice geometric realisation as
Kapovich--Millson bending flows \cite{KapovichMillson}.

\subsection*{Acknowledgements}
We thank A. M. Vershik for inspiring discussions at an early stage
of this project. We thank Pavel Etingof and Valerio Toledano Laredo
for useful comments and stimulating discussions. The second author
wishes to thank the Department of Mathematics of MIT for
hospitality. The third author would like to thank the Institute for
Mathematical Research at ETH Zurich for hospitality. We thank the
referee for providing relevant references and correcting misprints.
This work was partially supported by the Swiss National Science
Foundation (Grant 200020-122126).

\section{Classification of Gaudin subalgebras}
Since $\mathfrak t_n$ is defined by homogeneous relations, it is
graded in positive degrees: $\mathfrak t_n=\oplus_{i\geq1}\mathfrak
t_n^i$, with $\mathfrak t^1_n=\oplus_{i<j} \K t_{ij}$, $\mathfrak
t^2_n=\oplus_{i<j<k}\K [t_{ij},t_{ik}]$. In particular,
\[
\mathrm{dim}(\mathfrak t^1_n)=\frac{n(n-1)}2,\qquad
\mathrm{dim}(\mathfrak t^2_n)=\frac{n(n-1)(n-2)}6.
\]



\begin{defn}
A {\em Gaudin subalgebra} of $\mathfrak t_n$ is an abelian
subalgebra of maximal dimension contained in $\mathfrak t_n^1$.
\end{defn}
We will prove this maximal dimension is $n-1$. It follows from the
maximality condition that the central element
\[
c_n=\sum_{1\leq i<j\leq n} t_{ij}
\]
belongs to all Gaudin subalgebras.

\begin{exa}\label{ex-1}
The {\em Jucys--Murphy elements}
$t_{12},t_{12}+t_{13},\dots,\sum_{i=1}^{n-1} t_{in}$ are pairwise
commutative and thus span a Gaudin subalgebra. They play an
important role in the representation theory of the symmetric group,
see Remark \ref{r-1}.
\end{exa}
\begin{exa}
The main class of examples is provided by the spaces
\cites{Gaudin1976,Gaudin}
\begin{equation}\label{e-0}
G_n(z)=\left\{\sum_{1\leq i<j\leq
n}\frac{a_i-a_j}{z_i-z_j}\,t_{ij},\quad a\in \K^n\right\},
\end{equation}
parametrised by $z\in\Sigma_n/\mathrm{Aff}$, where
\[
\Sigma_n=\K^n\smallsetminus\cup_{i<j}\{z\in\K^n\,|\,z_i=z_j\}
\]
is the configuration space of $n$ distinct ordered points in the
plane and $\mathrm{Aff}$ is the group of affine maps $z\mapsto
az+b$, $a\neq 0$ acting diagonally on $\K^n$. This parameter space
is isomorphic to the moduli space $M_{0,n+1}= ((\mathbb
P^1)^{n+1}-\cup_{i<j}\{z_i=z_j\})/\mathit{PSL_2(\K)}$: the class of
$z$ is mapped to the class of $(z_1,\dots,z_n,\infty)$ in
$M_{0,n+1}$.
\end{exa}

\begin{lem}\label{l-1} The dimension of $G_n(z)$ is $n-1$.
\end{lem}

\begin{proof} The dimension is at most $n-1$ since there are $n$
parameters $a_1,\dots,a_n$ defined up to a common shift. Taking
$a=(1,\dots,1,0,\dots,0)$ with the number of ones ranging from 1 to
$n-1$ we obtain $n-1$ elements $K_j$ which are linearly independent:
$t_{j,j+1}$ appears in $K_j$ with non-vanishing coefficient but not
in $K_i$, $i\neq j$.
\end{proof}

The main result of this section is that Gaudin subalgebras are in
one to one correspondence with points in the Knudsen
compactification $\bar M_{0,n+1}$ of $M_{0,n+1}$, which is a
non-singular irreducible projective variety defined over $\mathbb Z$
\cite{Knudsen}. More precisely, we have the following result.

\begin{thm}\label{t-1}
Gaudin subalgebras in $\mathfrak t_n$ form a nonsingular subvariety
of the Grassmannian $G(n-1,n(n-1)/2)$ of $(n-1)$-planes in
$\mathfrak t_n^1$, isomorphic to $\bar M_{0,n+1}$.
\end{thm}

\begin{rem}\label{r-1}
To prove this theorem we only use the defining relations of
$\mathfrak t_n$ and the fact that both the generators $t_{ij}$,
$1\leq i<j\leq n$ and the brackets $[t_{ij},t_{ik}]$, $1\leq
i<j<k\leq n$ are linearly independent. Thus our result holds for any
quotient of $\mathfrak t_n$ with these properties. An important
example is the image of $\mathfrak t_n$ in the group algebra $\K
S_n$ of the symmetric group with commutator bracket, with $t_{ij}$
sent to the transposition of $i$ and $j$. An approach to the
representation theory of $S_n$ based on the simultaneous
diagonalization of the image of the Jucys--Murphy elements was
proposed in \cites{VershikOkounkov1996,VershikOkounkov}. Another
interesting case is given by the homomorphism $\phi\colon \mathfrak
t_n\mapsto U\mathfrak{o}(n)$ into the universal enveloping algebra
of the Lie algebra of the orthogonal group, sending $t_{ij}$ to
$X_{ij}^2$ where $X_{ij}$, $i<j$ are the standard generators of the
Lie algebra $\mathfrak o(n)$. The image consists of
(the complex versions of)
quantum Hamiltonians of the
corresponding Manakov tops \cite{Manakov}.
\end{rem}

The rest of this section is dedicated to the proof of Theorem
\ref{t-1}.

Let $D_n$ be the set of all distinct triples $(i,j,k)$ of numbers
between $1$ and $n$ (i.e., the set of injective maps
$\{1,2,3\}\to\{1,\dots,n\}$). For $(i,j,k)\in D_n$ denote by
$p_{ijk}\colon \mathfrak{t}_n^1\to \K^3$ the linear map
\[
\sum_{i<j}a_{ij}t_{ij}\to (a_{jk},a_{ik},a_{ij})
\]
where $a_{ij}$ is extended to all pairs by the rule $a_{ji}=a_{ij}$.
The map $p\colon D_n\to \mathrm{Hom}(\mathfrak{t}_n^1, \K^3), (i,j,k)\mapsto
p_{ijk}$ is equivariant under the natural action of $S_3$ on $D_n$
and on $\K^3$.

 We start with the following simple calculation, which was probably first done by Vinberg  \cite{Vinberg}.

\begin{lem}\label{l-2} Let $V\subset \mathfrak{t}_n^1$ be a Gaudin subalgebra.
Then, for all $(i,j,k)\in D_n$, $p_{ijk}(V)$ contains $(1,1,1)$ and
is at most two-dimensional.
\end{lem}

\begin{proof}
By the $S_3$-equivariance it is sufficient to prove the claim for
$i<j<k$. The space $p_{ijk}(V)$ contains $p_{ijk}(c_n)=(1,1,1)$. Let
$a=\sum_{i<j}a_{ij}t_{ij},b=\sum_{i<j}b_{ij}t_{ij}\in \mathfrak{t}_n^1$. The
commutator $[a,b]$ is a linear combination of the linearly
independent elements $[t_{ij},t_{jk}]$, $1\leq i<j<k\leq n$. Then
the equation $[a,b]=0$ is equivalent to the system
\[
  a_{ij}b_{jk}-a_{ij}b_{ik}
 +a_{ik}b_{ij}-a_{ik}b_{jk}
 +a_{jk}b_{ik}-a_{jk}b_{ij}=0,
\]
$1\leq i<j<k\leq n$. These equations are conveniently written in
determinant form
(cf. proof of Theorem 1 in \cite{Vinberg})
\begin{equation}\label{e-1}
\det\left(
\begin{array}{ccc}
 a_{jk}&b_{jk}&1\\
 a_{ik}&b_{ik}&1\\
 a_{ij}&b_{ij}&1
\end{array}
\right)=0.
\end{equation}
Thus $p_{ijk}(V)$ contains at most two linearly independent vectors.
\end{proof}

Thus for each Gaudin subalgebra there exist an $S_3$ equivariant map
$\ell\colon D_n\to (\K^3)^*$ sending $(i,j,k)$ to a linear form
$\ell_{ijk}$ vanishing on $(1,1,1)$ and such that
\begin{equation}\label{e-2}
\ell_{ijk}\circ p_{ijk}|_V=0.
\end{equation}
If $V=G_n(z)$ eq.~\eqref{e-2} is satisfied with the linear forms
\[
 \ell_{ijk}=(z_j-z_k,z_k-z_i,z_i-z_j).
\]
 Conversely, we have the following result.
\begin{lem}\label{l-3} Let $\ell\colon D_n\to (\K^3)^*$, $(i,j,k)\mapsto \ell_{ijk}$ be an
$S_3$-equivariant map such that $\ell_{ijk}(1,1,1)=0$ for all
$(i,j,k)$. Then
\[
V=\cap_{ijk}\mathrm{Ker}(\ell_{ijk}\circ p_{ijk}).
\]
is an abelian Lie subalgebra.
\end{lem}

\begin{proof} The vanishing condition implies that $c_n\in V$.
If $a,b\in V$ then $p_{ijk}(a),p_{ijk}(b)$ and $(1,1,1)$ belong to a
two-dimensional subspace of $\K^3$ and therefore obey \eqref{e-1}
for all $i,j,k$. It follows as in the proof of Lemma \ref{l-2} that
$[a,b]=0$.
\end{proof}

It remains to determine which systems of linear forms $\ell_{ijk}$
give commuting subspaces of maximal dimension. With respect to the
basis $t_{ij}$ of $\mathfrak{t}_n^1$ we can represent the linear forms
$\ell_{ijk}\circ p_{ijk}$ as the rows of a matrix $L$, so that the
corresponding commuting subspace is the kernel of $L$. The matrices
arising in this way belong to the following set.

\begin{defn}
Let $n\geq3$ and $\mathcal L_n$ be the space of matrices whose rows
are labeled by triples in $D^+_n=\{(i,j,k), 1\leq i<j<k\leq n\}$,
whose columns are labeled by pairs in $Z^+_n=\{(i,j), 1\leq i<j\leq
n\}$ and such that
\begin{enumerate}
\item The matrix elements in the row labeled by $(i,j,k)\in D^+_n$
vanish except possibly those in the columns $(j,k)$, $(i,k)$,
$(i,j)$.
\item Each row has at least a non-vanishing matrix element.
\item The sum of the matrix elements in each row is zero.
\end{enumerate}
\end{defn}
For example, matrices in $\mathcal L_4$ are of the form
\begin{equation}\label{e-3}
\begin{array}{cl}
\begin{array}{ccccccc}
 \phantom{123}
 &\phantom{a_{124}}
 &\phantom{a_{124}}
 &\phantom{a_{124}}
 &\phantom{a_{124}}
 &\phantom{a_{124}}
 &\phantom{a_{124}}\\
 &12&13&14&23&24&34
 \end{array}
 \\
\begin{array}{c}
123\\
124\\
134\\
234\end{array} \left(
\begin{array}{cccccc}
c_{123}&b_{123}&0     &a_{123}&0     &0     \\
c_{124}&0     &b_{124}&0     &a_{124}&0     \\
0     &c_{134}&b_{134}&0     &0     &a_{134}\\
0     &0     &0     &c_{234}&b_{234}&a_{234}
\end{array}
\right),
\end{array}
\end{equation}
with nonzero rows and zero row sums.

\begin{prop}\label{p-1}
Let $L\in\mathcal L_n$, $m\leq n$. Let $L_m$ be the matrix obtained
from $L$ by taking the matrix elements labeled by $D^+_m\times
Z^+_m\subset D^+_n\times Z^+_n$. Then $L_m\in\mathcal L_m$.
\end{prop}
\begin{proof} The claim is an easy consequence of the definition.
\end{proof}

\begin{lem}\label{l-5}
Let $L\in\mathcal L_n$ and $L_m$, $3\leq m\leq n$, the submatrix
with labels in $D^+_m\times Z^+_m$. Then
\[
 \mathrm{rank}(L)\geq
 \frac{(n-1)(n-2)}{2},\] with equality if and only if
\[
 \mathrm{rank}(L_m)=\frac{(m-1)(m-2)}{2},\qquad \text{for all
 $m=3,\dots,n$.}
\]
\end{lem}

\begin{proof}
We claim that there exists a row index set $I=I_3\cup
I_4\cup\dots\cup I_n\subset D^+_n$ such that
\begin{enumerate}
\item
For each $m$ the set $I_m$ has $m-2$ elements; they are of the form
$(i,j,m)$ for some $i<j<m$.
\item
For each $m$ there are distinct indices
$k_1,\dots,k_{m-1}\in\{1,\dots, m\}$ and an ordering $r_1,\dots,
r_{m-2}$ of $I_m$ such that the entry of row $r_i$ in column
$(k_j,m)$ is zero for $i<j$ and nonzero if $i=j$.
\end{enumerate}
The $(n-1)(n-2)/2$ rows of $L$ labeled by $I$ are then clearly
linearly independent, and the same holds for the rows of $L_m$ in
$I_3\cup\cdots\cup I_m$ for all $m\leq n$. It is also clear that if
a row of $L$ labeled by $D^+_m\subset D^+_n$ is a linear
combinations of rows labeled by $I$ then it is a linear combinations
of rows labeled by $I_m$. This proves the Lemma assuming the
existence of $I$.

To describe the construction of $I$ it is notationally convenient to
think of $D^+_n$ as the set $D_n/S_3$ of subsets of three elements
and thus to identify $(i,j,k)=(j,i,k)=(i,k,j)$. The row indices
$I_m$ can then be taken as $r_i =(\sigma_m(i),\sigma_m(i+1),m)$, for
some permutation $\sigma_m$ of $\{1,\dots,m-1\}$ such that the entry
of the row $r_j$ in the column $(\sigma_m(j+1),m)$ is nonzero. By
Property (1) of $\mathcal{L}_n$ this then implies that only $r_j$
among $r_1,\dots,r_{m-1}$ has a non zero entry in this column and we
set $k_j=\sigma_m(j+1)$, as desired. It remains to prove that such a
permutation exists. Let $\Gamma_m$ be the complete graph with vertex
set $\{1,\dots,m\}$. Pick an orientation $i\to j$ on each edge
$\{i,j\}$ of $\Gamma_m$ such that the entry in the row $(i,j,m)$ and
the column $(j,m)$ is nonzero. Such an orientation exists since at
least two of the entries in columns $(i,j)$, $(i,m)$, $(j,m)$ are
nonzero. Then the claim follows from the following simple  result of
elementary graph theory:
\begin{lem}\label{l-6}
For any orientation of the edges of a complete graph with $k$
vertices, there exists an oriented path
$\sigma(1)\to\sigma(2)\to\cdots\to\sigma(k)$ visiting each vertex
exactly once.
\end{lem}
In graph theory such a path is called Hamiltonian and this fact is
known as the existence of a Hamiltonian path in any tournament
\cite{Redei}.

The proof is by induction: for $k=1$ there is nothing to prove. Let
$\Gamma_k$ be the complete graph with vertex set $\{1,\dots,k\}$. An
orientation of its edges restricts to an orientation of the edges of
$\Gamma_{k-1}\subset \Gamma_k$. If we have a path $\gamma$ on
$\Gamma_{k-1}$ starting at a vertex $i$ and ending at a vertex $j$
then either there is an edge $k\to i$ or $j\to k$ and we can
complete $\gamma$ to a path in $\Gamma_k$ by adding it, or there
exists a step $a\to b$ of $\gamma$ that can be replaced by $a\to
k\to b$ to obtain a path in $\Gamma_k$ with the required property.
\end{proof}

\begin{cor}\label{c-1} Abelian subalgebras lying in $\mathfrak{t}_n^1$ have dimensions at most
$n-1$.
\end{cor}

Indeed the rank of a matrix in $\mathcal L_n$ is at least
$(n-1)(n-2)/2$ and its kernel has dimension at most
\[
\frac{n(n-1)}{2}-\frac{(n-1)(n-2)}{2}=n-1.
\]

\begin{lem}\label{l-7}
Suppose $L\in\mathcal L_n$ has minimal rank $(n-1)(n-2)/2$, so that
$\mathrm{Ker}\,L$ defines a Gaudin subalgebra. Denote the entries of
row $(i,j,k)$ in columns $(i,j),(i,k),(j,k)$ by
$a_{ijk},b_{ijk},c_{ijk}$ respectively.  Then
\begin{gather}
\label{e-3.5}
(a_{ijk}:b_{ijk}:c_{ijk})=(b_{jik}:a_{jik}:c_{jik})=(a_{ikj}:c_{ikj}:b_{ikj})
\\
 \label{e-4}
a_{ijk}+b_{ijk}+c_{ijk}=0,
\end{gather}
for all $(i,j,k)\in D_n$, and
\begin{equation}\label{e-5}
b_{ijk}c_{ijl}b_{ikl}+c_{ijk}b_{ijl}c_{ikl}=0,
\end{equation}
for all $(i,j,k,l)\in V_n$.
\end{lem}
\begin{proof} The first two equations are a rephrasing of the
$S_3$-equivariance and the condition $\ell_{ijk}(1,1,1)=0$. Consider
the third equation. By possibly renumbering the vertices we can
assume that the four indices are 1,2,3,4. By Lemma \ref{l-5} the
submatrix $L_4$ with labels in $D^+(4)\times Z^+(4)$ has rank $3$.
This matrix has the form \eqref{e-3}. Since the last row is nonzero,
the upper left $3\times 3$ minor vanishes.
\end{proof}

We can now conclude the proof of Theorem \ref{t-1}. Let $Z_n$ be the
subvariety of Gaudin subalgebras in the Grassmannian of $n-1$ planes
in $\mathfrak{t}_n^1$. By Lemma \ref{l-2}, every Gaudin subalgebra
$V$ is contained in $\mathrm{Ker}\,L$ for some $L\in\mathcal L_n$.
Since by Lemma \ref{l-3} $\mathrm{Ker}\,L$ is an abelian subalgebra,
it follows by maximality that $V$ is actually equal to
$\mathrm{Ker}\,L$. Let $Y_n$ be the subvariety of $(\mathbb
P^2)^{D_n}$ defined by the equations \eqref{e-3.5}--\eqref{e-5} for
homogeneous coordinates $(a_{ijk}:b_{ijk}:c_{ijk})$. By Lemma
\ref{l-7} we then have a map $Z_n\to Y_n$ which is clearly
injective. Since the subalgebras $G_n(z)$ are contained in $Z_n$,
$Z_n$ has a component of dimension $n-2$. As we will presently see,
$Y_n$ is isomorphic to the non-singular irreducible projective
$(n-2)$-dimensional variety $\bar M_{0,n+1}$ and therefore $Z_n\to
Y_n$ is an isomorphism.

It remains to prove that $Y_n\simeq \bar M_{0,n+1}$. In order to do
so let us first notice that by \eqref{e-4}, $Y_n$ actually lies in
$(\mathbb P^1)^{D_n}$, where $\mathbb P^1\subset\mathbb P^2$ is
embedded as $(x:y)\mapsto (x-y:-x:y)$. It is easy to rewrite the
remaining equations defining $Y_n$ in these coordinates:
\begin{enumerate}
\item
$x_{ikj}x_{ijk}=y_{ikj}y_{ijk}$,
\item
$x_{jik}y_{ijk}=y_{ijk}y_{jik}-x_{ijk}y_{jik}$,
\item
$x_{ijk}y_{ijl}x_{ikl}=y_{ijk}x_{ijl}y_{ikl}$.
\end{enumerate}
It turns out that, by a variant of a Theorem of Gerritzen, Herrlich
and van der Put, these are precisely the relations defining $\bar
M_{0,n+1}$ as a subvariety of $(\mathbb P^1)^{D_n}$. We deduce this
variant from the original Theorem of Gerritzen, Herrlich and van der
Put in the Appendix, see Theorem \ref{t-3}.

\section{The sheaf of Gaudin subalgebras}
By Theorem \ref{t-1} Gaudin subalgebras of $\mathfrak t_n$ form a
family of vector spaces $\mathcal G_n$ on $\bar M_{0,n+1}$: the
fibre at $z$ is the Gaudin subalgebra corresponding to $z$. The
purpose of this section is to identify this family in terms of the
geometry of $\bar M_{0,n+1}$.

Consider first the Gaudin subalgebras parametrised by $M_{0,n+1}$.
Let $\tilde M_{0,n+1}=\Sigma_n/\K$, where the group $\K\subset
\mathrm{Aff}$ is the translation subgroup. The natural projection
$p\colon \tilde M_{0,n+1}\to M_{0,n+1}$ is a principal
$\Kstar$-bundle. The $\mathfrak t^1_n$-valued 1-form
\[
\omega=\sum_{i<j}\frac {dz_i-dz_j}{z_i-z_j}t_{ij},
\]
is a $\Kstar$-invariant element $\Omega^1(\tilde
M_{0,n+1})\otimes\mathfrak t^1_n$. The pairing with $\omega$ defines
a map
\[
T_z{\tilde M_{0,n+1}}\to \mathfrak t_n^1
\]
from the tangent space at $z\in \tilde M_{0,n+1}$ to $\mathfrak
t_n^1$ with image $G_{n}(z)$, see \eqref{e-0}. By Lemma \ref{l-1}
this map is injective. Now $G_n(z)=G_n(z')$  if and only if
$z'=\lambda z$ for some $\lambda\in\Kstar$. More precisely the
action of $\Kstar$ on $\tilde M_{0,n+1}$ lifts naturally to
$\mathit{T\tilde{M}}_{0,n+1}$ and the invariance of $\omega$ implies
that $\omega$ defines an injective bundle map
\[
\begin{array}{ccc}
\mathit{T\tilde{M}}_{0,n+1}/\Kstar&\to &M_{0,n+1}\times \mathfrak
t^1_n\\
\downarrow&&\downarrow\\
M_{0,n+1}&=&M_{0,n+1}
\end{array}
\]
Its image is a vector bundle with fibres $G_n(z)$, $z\in M_{0,n+1}$.
Moreover  $\mathit{T\tilde M}_{0,n+1}/\Kstar$ is an extension of the
tangent bundle to $M_{0,n+1}$: the kernel of the natural surjective
bundle map $\mathit{T\tilde M}_{0,n+1}/\Kstar\to
\mathit{TM}_{0,n+1}$ is spanned by the class of the Euler vector
field
\[
E=\sum_{i=1}^n z_i\frac\partial{\partial z_i},
\]
generating the $\Kstar$-action. Turning to the language of sheaves,
more convenient when we pass to $\bar M_{0,n+1}$, we thus have an
exact sequence of locally free sheaves on $X=M_{0,n+1}$
\begin{equation}\label{e-101}
0\to \mathcal O_X\to \mathcal G_n\to T_X\to 0.
\end{equation}
Here $T_X$ is the sheaf of vector fields and $\mathcal G_n$ is the
sheaf of $t^1_n$-valued functions whose value at each $z$ lies in
$G_n(z)$. For any open set $U\subset M_{0,n+1}$, $\mathcal G_n(U)$
may be identified via the map $\omega$ with the space of
$\Kstar$-invariant vector fields on $p^{-1}(U)$.

As is well-known, invariant vector fields can be identified with
first order twisted differential operators on the base manifold:

\begin{lem}\label{l-twisted}
Let $p\colon P\to X$ be a principal $\Kstar$-bundle on a smooth
variety $X$, $L=P\times_{\Kstar}\K$ be the associated line bundle
where $\Kstar$ acts on $\K$ by multiplication. Then for each open
set $U\subset X$, the Lie algebra $T_P(p^{-1}(U))^{\Kstar}$ is
isomorphic to the Lie algebra $D^1_{L^\vee}(U)$ of first order
differential operators acting on sections of the dual line bundle
$L^\vee$ (i.e., twisted by $L^\vee$).
\end{lem}

\begin{proof}
A section of $L^\vee$ on $U$ is the same as a function $f\colon
p^{-1}(U)\to \K$ such that $f(y\cdot\lambda)=\lambda f(y)$, $y\in
p^{-1}(U)$, $\lambda\in \Kstar$. It is clear from this
representation that $T_P(p^{-1}(U))^{\Kstar}$ acts on sections of
$L^\vee$. Moreover the infinitesimal generator $E$ of the
$\Kstar$-action acts by 1. Thus, upon choosing a local
trivialization of $P$, we may write any invariant vector field as
$\xi+f E$ where $\xi$ is a vector field on $U$ and $f\in O_X(U)$.
This invariant vector field acts on a section as the first order
differential operator $\xi+f$.
\end{proof}

As a consequence we have a description of $\mathcal G_n$ as a sheaf
of twisted first order differential operators:

\begin{prop}
Gaudin subalgebras corresponding to points of $M_{0,n+1}$ form a
locally free sheaf isomorphic to the sheaf $D^1_{L^\vee}$ of first
order differential operators on $M_{0,n+1}$ twisted by the line
bundle $L^\vee$ dual to the associated bundle to $\tilde M_{0,n+1}$
via the identity character of $\Kstar$.
\end{prop}

Let us now extend this to $\bar M_{0,n+1}$.

Let us first recall some known facts about the geometry of $\bar
M_{0,n+1}$ \cites{Knudsen,Gerritzenetal, Keel}. This space can be
defined as the closure in $(\mathbb P^1)^{D_n}$ of the image of the
injective map
\[
\mu\colon M_{0,n+1}=\Sigma_n/\mathrm{Aff}\to(\mathbb P^1)^{D_n},
\]
sending the class of $z\in\Sigma_n$ to the collection of cross
ratios involving the point at infinity
\[
\mu_{ijk}(z)=\frac{z_i-z_k}{z_i-z_j}=\frac{(z_i-z_k)(\infty-z_j)}{(z_i-z_j)(\infty-z_k)},\qquad
(i,j,k)\in D_n.
\]
Moreover the image is characterized by an explicit set of equations,
see Theorem \ref{t-3}.

The complement of $M_{0,n+1}$ in $\bar M_{0,n+1}$ is a normal
crossing divisor $D=\cup D_S$, where the union is over all subsets
$S$ of $\{1,\dots,n\}$ with at least two and at most $n-1$ elements.
The irreducible component $D_S$ is isomorphic to $\bar
M_{0,m+1}\times \bar M_{0,n-m+2}$ and is the closure of the
subvariety consisting of stable curves with one nodal point such
that the marked points labeled by $S$ are those on the component not
containing the point labeled by $n+1$.

Local coordinates on a neighbourhood in $\bar M_{0,n+1}$ of a
generic point of $D_S$ are the cross ratios $\zeta_r=\mu_{ijr}$,
$r\in S\smallsetminus\{i,j\}$, $z_s=\mu_{iks}$, $s\in
S^c\smallsetminus\{k\}$, $t=\mu_{ijk}$, for any fixed $i,j\in S$,
$k\not\in S$ ($S^c$ denotes the complement of $S$ in
$\{1,\dots,n\}$). In these coordinates, $D_S$ is given by the
equation $t=0$. The change of variables from these coordinates to
the coordinates $z_i$ of $M_{0,n+1}$ is as follows. We may assume
that $S=\{1,\dots,m\}$ and choose $i=1,j=m,k=m+1$ (the general case
can be obtained from this by permuting the coordinates). Then for
generic $t\neq0$ the point in $M_{0,n+1}$ with coordinates
$(\zeta,z,t)$ is
\begin{equation}\label{e-6}
[0,t\zeta_2,\dots,t\zeta_{m-1},t,z_2,\dots,z_{n-m},1]\in
M_{0,n+1}=\Sigma_n/\mathrm{Aff}.
\end{equation}

To extend the exact sequence \eqref{e-101} to $\bar M_{0,n+1}$ we
first show that $\tilde M_{0,n+1}$ extends to a principal
$\Kstar$-bundle
\[
p\colon\tilde{\bar M}_{0,n+1}\to \bar M_{0,n+1}.
\]
This can be seen from the presentation of Theorem \ref{t-3}. Let $H$
be the kernel of the product map $(\Kstar)^{D_n}\to \Kstar$. Then
$P_n=(\K^2\smallsetminus\{0\})^{D_n}/H$ is a principal
$\Kstar$-bundle on $(\mathbb P^1)^{D_n}$ and $\tilde M_{0,n+1}$
embeds into $P_n$ (via $z\mapsto \mbox{class of
}((z_i-z_k,z_i-z_j)_{(i,j,k)\in D_n})$) as the restriction of $P_n$
to the image of $M_{0,n+1}$ in $(\mathbb P^1)^{D_n}$. We then define
$\tilde{\bar M}_{0,n+1}$ to be the restriction of $P_n$ to $\bar
M_{0,n+1}\subset (\mathbb P^1)^{D_n}$.

Recall that the locally free sheaf $T_X\langle-D\rangle$ of
logarithmic vector fields on a variety $X$ with a normal crossing
divisor $D$ consists of vector fields whose restriction to a generic
point of $D$ is tangent to $D$. It is dual to the sheaf
$\Omega^1_X\langle D\rangle$ of logarithmic 1-forms, spanned over
$\mathcal O_X$ by regular 1-forms and $df/f$ where $f\in \mathcal
O_X$ with $f\neq 0$ on $X\smallsetminus D$, see
\cite{Deligne}*{Sect.~II.3}.

\begin{thm}\label{t-102}
Let $L$ be the associated line bundle $\tilde {\bar
M}_{0,n+1}\times_{\Kstar}\K$ with the identity character of
$\Kstar$. Gaudin subalgebras form a vector bundle on $\bar
M_{0,n+1}$. As a locally free sheaf, it is isomorphic to the sheaf
$D^1_{L^\vee}\langle -D\rangle$ of first order differential
operators on $\bar M_{0,n+1}$ twisted by $L^\vee$, whose symbol is
logarithmic. In particular there is an exact sequence of sheaves on
$X=\bar M_{0,n+1}$
\[
0\to \mathcal O_X\to \mathcal G_n\to T_X\langle -D\rangle\to 0.
\]
The embedding of the trivial bundle $\mathcal O_X$ sends $1$ to
$c_n=\sum_{i<j}t_{ij}$.
\end{thm}

\begin{proof}
Let us introduce the abbreviated notation $X=\bar M_{0,n+1}$,
$p\colon \tilde X=\tilde{\bar M}_{0,n+1}\to X$. Let $\tilde
D=p^{-1}(D)$ be the pull-back to $\tilde X$ of the divisor $D$. Then
$\omega$ is a form in $\Omega^1(\tilde X)\otimes \mathfrak t^1_n$
with logarithmic coefficients. Indeed, in the coordinates
$\zeta,z,t$ of \eqref{e-6} and fibre coordinate $\lambda\in\Kstar$
around a generic point of $\tilde D$, $\omega$ has the local
coordinate expression
\begin{equation}\label{e-105}
\begin{array}{rcl}
 \omega&=&\sum_{1\leq i<j\leq m} t_{ij}d\log(\zeta_i-\zeta_j)
 +\sum_{1\leq i<j\leq m}t_{ij}d\log(t)
 \\
 &&
 \\
 &&
 +\sum_{m<i<j\leq n}t_{ij} d\log(z_i-z_j)
 +\sum_{1\leq i<m<j\leq n}t_{ij}d\log(z_j)
 \\
 &&
 \\
 &&
 +\,c_n d\log(\lambda)
\end{array}
\end{equation}
with the understanding that $\zeta_0=0, \zeta_m=z_n=1$. Thus
$\omega$ may be paired with invariant logarithmic vector fields on
$\tilde X$, which in turn may be identified by Lemma \ref{l-twisted}
with first order differential operators on $X$, to give a map $
D^1_{L^\vee}\langle -D\rangle\to \mathfrak t^1\otimes \mathcal O_X,
$ which is injective on $M_{0,n+1}$. We need to show that the map is
injective on all of $X=\bar M_{0,n+1}$. Since the locus of
non-injectivity is empty or of codimension one, it is sufficient to
show that the map is injective as we approach a generic point of the
divisor $D$. This is easy to check using \eqref{e-105}.

The embedding of $\mathcal O_X$ sends $1$ to
$\langle\omega,E\rangle=c_n$.
\end{proof}

Thus $\mathcal G_n$ is a sheaf of twisted first order differential operators with
regular singularities along $D$, see \cite{Ueno}, Sect.~5.2.

\begin{rem} The divisor class of the line bundle $L$ can be easily computed by choosing
a section. The result is
\[
[L]=-\sum_{S\supset\{1,2\}}[D_S].
\]
\end{rem}

\section{Coisotropic spectra} Suppose that $U\mathfrak t_n\to A$ is
a homomorphism of unital algebras, e.g. $A=\mathrm{End}(V)$ for some
representation $V$ of $\mathfrak t_n$. Then we get a sheaf $\mathcal
E_A$ of commutative subalgebras of $A$ on $X=\bar M_{0,n+1}$ as the
image of the symmetric $\mathcal O_X$-algebra $S\mathcal G_n$. By
Theorem \ref{t-102}, $\mathcal G_n$ is naturally a sheaf of Lie
algebras, so that $S\mathcal G_n$ is a sheaf of Poisson algebras.

\begin{cor}\label{c-103} Let $\varphi\colon U\mathfrak t^1_n\to A$ be an algebra homomorphism. Then
the kernel of the induced map of sheaves of algebras $ \varphi\colon
S\mathcal G_n\to A\otimes \mathcal O_X, $ is closed under Poisson
brackets.
\end{cor}

\begin{proof} This is basically a consequence of the fact that
$\omega$ is closed. The map is defined by identifying $\mathcal G_n$
with the sheaf of $\Kstar$-invariants of $p_*T_{\tilde X}\langle
-\tilde D\rangle$. It is the algebra homomorphism sending an
invariant section $\xi$ to  $\omega_A(\xi)=\varphi(\omega(\xi))$. A
Poisson bracket of monomials $\xi_1\cdots\xi_k$,
$\eta_1,\cdots,\eta_l$, $\xi_i,\eta_j\in \mathcal G_n\simeq
p_*T_{\tilde X}\langle -\tilde D\rangle^{\Kstar}$ is sent to
\begin{eqnarray*}
\lefteqn{\varphi(\{\xi_1\cdots\xi_k,\eta_1\cdots\eta_l\})=
 \sum_{i,j}\varphi([\xi_i,\eta_j]\xi_1\cdots\hat\xi_i\cdots\xi_k\,\eta_1\cdots\hat\eta_j\cdots\eta_l)
 }
 \\
 &=&
 \sum_{i,j}\omega_A([\xi_i,\eta_j])\prod_{r\neq i}\omega_A(\xi_r)\prod_{s\neq j}\omega_A(\eta_s)
 \\
 &=&
 \sum_{i,j}\prod_{r<i}\omega_A(\xi_r)\prod_{s<j}\omega_A(\eta_s)
 \,\omega_A([\xi_i,\eta_j])\prod_{r>i}\omega_A(\xi_r)\prod_{s>j}\omega_A(\eta_s)
 \\
 &=&
 \sum_{i,j}\prod_{r<i}\omega_A(\xi_r)\prod_{s<j}\omega_A(\eta_s)
 \,(\xi_i\omega_A(\eta_j)-\eta_j\omega_A(\xi_i))\prod_{r>i}\omega_A(\xi_r)\prod_{s>j}\omega_A(\eta_s)
 \\
 &=&
 \sum_{i,j}\prod_{r<i}\omega_A(\xi_r)\prod_{s<j}\omega_A(\eta_s)
 \,\xi_i\omega_A(\eta_j)\prod_{s>j}\omega_A(\eta_s)\prod_{r>i}\omega_A(\xi_r)
 \\
 & &-
 \sum_{i,j}\prod_{s<j}\omega_A(\eta_s)\prod_{r<i}\omega_A(\xi_r)
 \,\eta_j\omega_A(\xi_i)\prod_{r>i}\omega_A(\xi_r)\prod_{s>j}\omega_A(\eta_s)
 \\
 &=&
 \sum_{i}\prod_{r<i}\omega_A(\xi_r)\,\xi_i\varphi(\eta_1\cdots\eta_l)
 \prod_{r>i}\omega_A(\xi_r)
 \\
 & &
 -
 \sum_{j}\prod_{s<j}\omega_A(\eta_s)
 \eta_j\varphi(\xi_1\cdots\xi_k)\prod_{s>j}\omega_A(\eta_s).
\end{eqnarray*}
In this calculation we use the fact that
$\omega_A(\xi_r),\omega_A(\eta_s)$ commute with each other; the
peculiar choice of ordering of factors is necessary since their
derivatives $\xi_i\omega_A(\eta_j)$, $\eta_j\omega_A(\xi_i)$ do not
necessarily commute with them.

It follows that if $\varphi(a)=\varphi(b)=0$ then also
$\varphi(\{a,b\})=0$, which is the claim.
\end{proof}

Recall that $\Omega^1_{\tilde X}\langle \tilde D\rangle$ is locally
free and thus the sheaf of sections of a vector bundle, the
logarithmic cotangent bundle. Let us denote by $T^*\tilde X\langle
\tilde D\rangle$ the total space of this vector bundle.  It is the
relative spectrum of the symmetric algebra $ST_{\tilde
X}\langle-\tilde D\rangle$, which is a sheaf of Poisson algebras.
Thus $T^*\tilde X\langle \tilde D\rangle$ is a Poisson variety; the
Poisson structure restricts to the usual symplectic structure on the
cotangent bundle of $\tilde M_{0,n+1}$. The group $\Kstar$ acts on
$\tilde X$ preserving $\tilde D$. This action lifts to a Hamiltonian
action on the logarithmic cotangent bundle with moment map
$E\in\Gamma(\tilde X,T_{\tilde X}\langle -\tilde D\rangle)\subset
\Gamma(\tilde X,ST_{\tilde X}\langle -\tilde D\rangle)=\mathcal
O(T^*\tilde X\langle \tilde D\rangle)$.

\begin{defn} Let $\alpha\in\K$. The {\em twisted logarithmic
cotangent bundle} $T^*X_\alpha\langle D\rangle$ {\em with twist
$\alpha$} is the Hamiltonian reduction $E^{-1}(\alpha)/\Kstar$.
\end{defn}

By construction $T^*X_\alpha\langle D\rangle$ is a Poisson variety.
For $\alpha=0$ it is the logarithmic cotangent bundle of $X$. By
definition, the regular functions on an open set $\tilde
U=p^{-1}(U)$, $U\subset X$, are section of $(S T_{\tilde X}\langle-
D\rangle(\tilde U))^{\Kstar}/I(\tilde U)$ where $I(\tilde U)$ is the
ideal generated by $E-\alpha$.

Then Corollary \ref{c-103} can be reformulated as follows.

\begin{cor}\label{c-104}
Let $U\mathfrak t^1_n\to A$ be an algebra homomorphism such that
$c_n$ is mapped to $\alpha 1$ and let $\mathcal E_A$ be the
corresponding sheaf of commutative algebras on $X=\bar M_{0,n+1}$.
Then the relative spectrum of $\mathcal E_A$ is a coisotropic
subscheme of the twisted logarithm cotangent bundle $T^*X_\alpha
\langle D\rangle$.
\end{cor}

In particular, if $A$ is finite dimensional, then the part of the
spectrum over $X^0=M_{0,n+1}$ is a Lagrangian subvariety of the
symplectic manifold $T^*X^0_\alpha=(E^{-1}(\alpha)\cap T^*\tilde
X^0)/\Kstar$.

\appendix

\section{Stable curves of genus zero}
Recall that a stable curve of genus zero with $r\geq 3$ marked
points is a pair $(C,S)$ where $C$ is a connected projective
algebraic curve of genus $0$ whose singularities are ordinary double
points and $S=(p_1,\dots,p_r)$ is an ordered set of distinct
nonsingular points of $C$ such that each irreducible component of
$C$ has at least three special (marked or singular) points. The
genus zero condition means that the irreducible components are
projective lines whose intersection graph is a tree. The moduli
space $\bar M_{0,r}$ of stable rational curves with $r\geq 3$ marked
point \cite{Knudsen} is a smooth algebraic variety of dimension
$r-3$ defined over $\mathbb Q$. It contains as a dense open set the
quotient of the configuration space
\[
 M_{0,r}=\{z\in (\mathbb P^1)^r\,|\, z_i\neq z_j,
 \text{for all $i\neq j$}\}/\mathit{PSL}_2.
\]
of $r$ distinct labeled points on the projective line by the
diagonal action of $\mathrm {Aut}(\mathbb P^1)\simeq
\mathit{PSL}_2$.

Here is a simple description of $\bar M_{0,r}$, due to Gerritzen,
Herrlich and van der Put \cite{Gerritzenetal}. For each distinct
$(i,j,k,l)$ in $\{1,\dots,r\}$, let $\lambda_{ijkl}\colon M_{0,r}\to
\mathbb P^1$ be the map sending the class of $z$ to the cross-ratio
 \[
 \lambda_{ijkl}(z)=\frac{(z_i-z_l)(z_j-z_k)}{(z_i-z_k)(z_j-z_l)}
 \in \K\subset \mathbb P^1.
 \]
Let $V_r$ be the set of distinct quadruples of integers between 1
and $r$. Then $\bar M_{0,r}$ is the closure of the image of the
embedding
\[
M_{0,r}\to (\mathbb P^1)^{V_r},\quad z\to (\lambda_v(z))_{v\in V_r},
\]
sending $z$ to the system of cross-ratios
$\lambda_{ijkl}(z)\in\mathbb P^1\smallsetminus\{0,\infty,1\}$.

The cross-ratios $\lambda_v=\lambda_v(z)$ of a point in $M_{0,r}$
obey the relations
\begin{enumerate}
\item[($\lambda1$)] $\lambda_{jikl}=1/\lambda_{ijkl}$ for all distinct $i,j,k,l$.
\item[($\lambda2$)] $\lambda_{jkli}=1-\lambda_{ijkl}$ for all distinct $i,j,k,l$.
\item[($\lambda3$)]  $\lambda_{ijkl}\lambda_{ijlm}=\lambda_{ijkm}$ for all distinct
$i,j,k,l,m$.
\end{enumerate}

\begin{thm}\label{t-2}(Gerritzen, Herrlich, van der Put \cite{Gerritzenetal})
The subvariety of $(\mathbb P^1)^{V_r}$ defined by these relations,
more precisely by their version for homogeneous coordinates
$\lambda_{v}=x_v/y_v$:
\begin{enumerate}
\item
$x_{jikl}x_{ijkl}=y_{jikl}y_{ijkl}$ for all distinct $i,j,k,l$.
\item
$x_{jkli}y_{ijkl}=y_{ijkl}y_{jkli}-x_{ijkl}y_{jkli}$ for all
distinct $i,j,k,l$.
\item
$x_{ijkl}x_{ijlm}y_{ijkm}=y_{ijkl}y_{ijlm}x_{ijkm}$ for all distinct
$i,j,k,l,m$,
\end{enumerate}
is a fine moduli space of stable curves of genus zero. The dense
open subvariety $M_{0,r}$ is embedded via the cross-ratios $z\mapsto
(x_v:y_v)_{v\in V_r}$, with
\[
(x_{ijkl}:y_{ijkl})=((z_i-z_l)(z_j-z_k):(z_i-z_k)(z_j-z_l)).
\]
\end{thm}

For our purpose it is useful to have a more economical description
of the moduli space by taking only cross-ratios involving a
distinguished marked point. Let $n=r+1\geq2$ and $D_n$ be the set of
all distinct triples $(i,j,k)$ of integers between $1$ and
$n$.\footnote{Following \cite{Gerritzenetal}, we denote the sets of
distinct pairs, triples and quadruples by the initials $Z$, $D$, $V$
of the corresponding German or Dutch numerals}\label{f-1} The
cross-ratios $\mu_{ijk}(z)=\lambda_{i,n+1,j,k}(z)$ obey
\begin{enumerate}
\item[($\mu1$)]
$\mu_{ikj}=1/\mu_{ijk}$, for all distinct $i,j,k$,
\item[($\mu2$)]
$\mu_{jik}=1-\mu_{ijk}$,  for all distinct $i,j,k$,
\item[($\mu3$)]
$\mu_{ijk}\mu_{ikl}=\mu_{ijl}$, for all distinct $i,j,k,l$.
\end{enumerate}
The claim is that the homogeneous version of these relations define
$\bar M_{0,n+1}$:
\begin{thm}\label{t-3}
The moduli space $\bar M_{0,n+1}$ is isomorphic to the subvariety of
$(\mathbb P^1)^{D_n}$ defined by the equations
\begin{enumerate}
\item
$x_{ikj}x_{ijk}=y_{ikj}y_{ijk}$ for all $(i,j,k)\in D_n$.
\item
$x_{jik}y_{ijk}=y_{ijk}y_{jik}-x_{ijk}y_{jik}$,  for all $(i,j,k)\in
D_n$,
\item
$x_{ijk}x_{ikl}y_{ijl}=y_{ijk}y_{ikl}x_{ijl}$ for all $(i,j,k,l)\in
V_n$.
\end{enumerate}
The open subvariety $M_{0,n+1}$ is embedded via the cross-ratios
$z\mapsto (x_d:y_d)_{d\in D_n}$, with
\[
(x_{ijk}:y_{ijk})=((z_i-z_k)(z_{n+1}-z_j):(z_i-z_j)(z_{n+1}-z_k))
\]
\end{thm}

\begin{rem} In \cite{Gerritzenetal}, $\bar M_{0,n+1}$ is considered as a
scheme over $\mathbb Z$, being defined as a subscheme of
$\prod_{v\in V_{n+1}}\mathrm{Proj}(\mathbb Z[x_v,y_v])$. Our proof
applies also to this more general setting. \end{rem}

\begin{proof}
Let us denote by $Y_n$ the subvariety of $(\mathbb P^1)^{D_n}$
defined by these relations. We have an obvious map $f\colon\bar
M_{0,n+1}\to Y_n$, the projection onto the cross-ratios
$\mu_{ijk}=(x_{i,n+1,j,k}:y_{i,n+1,j,k})$ with $(i,j,k)\in D_n$. We
show that this map is an isomorphism by constructing the inverse map
$\mu\mapsto \lambda=g(\mu)$. If one of $i,j,k,l$ is equal to $n+1$,
$\lambda_{ijkl}$ is obtained using ($\lambda1$) and ($\lambda2$)
from
\begin{equation}\label{e-A1}
\lambda_{i,n+1,k,l}=\mu_{ikl}
\end{equation}
For $(i,j,k,l)\in V_n$,
$\lambda_{ijkl}$ is given by either
\begin{equation}\label{e-A2}
\lambda_{ijkl}=\frac{\mu_{ikl}}{\mu_{jkl}}
\end{equation}
 or
\begin{equation}\label{e-A3}
\lambda_{ijkl}=\frac{\mu_{kij}}{\mu_{lij}},
\end{equation}
depending on which of the two expressions is defined (i.e. not $0/0$
or $\infty/\infty$.)\footnote{The notation $x_3=x_1/x_2$ for
$x_i=(x'_i:x''_i)\in\mathbb P^1$ means
$x'_3x'_2x''_1=x_1'x''_2x''_3$; this defines $x_3$ given $x_1$ and
$x_2$ unless $x_1$ and $x_2$ are both  $0=(0:1)$ or $\infty=(1:0)$}

We first check that \eqref{e-A1}--\eqref{e-A3} correctly define a
map $Y_n\to (\mathbb P^1)^{V_{n+1}}$. First of all $(\lambda1)$ and
$(\lambda2)$ say that the map $\lambda\colon V_{n+1}\to \mathbb P^1$
is equivariant under the natural action of $S_4$ on $V_{n+1}$ and on
$\mathbb P^1$ by fractional linear transformations. Since, by
$(\mu1)$ and $(\mu2)$, $\mu$ is $S_3$-equivariant, \eqref{e-A1}
defines consistently $\lambda_v$ for $v$ in the $S_4$-orbit of
$(i,n+1,k,l)$. Next, we claim that at least one of the ratios in
\eqref{e-A2}, \eqref{e-A3} is defined. Indeed suppose that
\eqref{e-A2} is not defined because $\mu_{ikl}=\mu_{jkl}=0$. Then,
by ($\mu1$) and ($\mu2$), $\mu_{kil}=1$ and
$\mu_{klj}=\mu_{kjl}^{-1}=(1-\mu_{jkl})^{-1}=1$. Thus, by $(\mu3)$,
$\mu_{kij}=1$ and \eqref{e-A3} is defined. Similarly, if
$\mu_{ikl}=\mu_{jkl}=\infty$ then $\mu_{lij}=1$. The same arguments
shows that $\eqref{e-A2}$ is defined if $\eqref{e-A3}$ is not. It
remains to show that if the right-hand sides of both \eqref{e-A2},
\eqref{e-A3} are defined then they are equal. The following identity
is useful for this purpose.

\begin{lem}\label{l-A1}
Suppose $(\mu_d)_{d\in D_n}$ obey $(\mu1)$, $(\mu2)$ and let
$(i,j,k)\in D_n$. Then either $\mu_{ijk},\mu_{jki},\mu_{kij}$ are a
permutation of $0,\infty,1$ or none of them belongs to
$\{0,\infty,1\}$. In the latter case their product is $-1$ or,
equivalently,
\[
\mu_{ijk}=-\frac{\mu_{jik}}{\mu_{kij}}.
\]
\end{lem}
Indeed, if $\mu_{ijk}=x$ then $\mu_{jki}=1/\mu_{jik}=1/(1-x)$ and
$\mu_{kij}=1-1/x$, which implies the Lemma.

This Lemma combined with $(\mu1), (\mu3)$ implies the following two
identities (holding whenever the expressions are defined):
\begin{equation}\label{e-A4}
\frac{\mu_{ikl}}{\mu_{jkl}}=\frac{\mu_{ikj}\,\mu_{ijl}}{\mu_{jki}\,\mu_{jil}}=
\frac{\mu_{kij}}{\mu_{lij}} ,\qquad
\frac{\mu_{kij}}{\mu_{lij}}=\frac{\mu_{kil}\,\mu_{klj}}{\mu_{lik}\,\mu_{lkj}}=
\frac{\mu_{ikl}}{\mu_{jkl}}.
\end{equation}
Assume that both right-hand sides of \eqref{e-A2},\eqref{e-A3} are
defined. We have four cases: (a)
$\mu_{ikl},\mu_{jkl}\notin\{0,\infty,1\}$. Then the second identity
in \eqref{e-A4} proves that \eqref{e-A2},\eqref{e-A3} agree. (b)
$\mu_{kij},\mu_{lij}\notin\{0,\infty,1\}$. Here the first identity
implies the claim. (c) $\mu_{ikl}=0$. Since, by $(\mu3)$,
$\mu_{ikl}=\mu_{ikj}\mu_{ijl}$, we have either $\mu_{ikj}=0$ or
$\mu_{ijl}=0$. In the first case $\mu_{kij}=1$ and therefore
$\mu_{kjl}=\mu_{kji}\mu_{kil}=1\cdot1=1$, implying $\mu_{jkl}=0$, in
contradiction with the assumption that the right-hand side of
\eqref{e-A2} is defined. In the second case,
$\mu_{lij}=1-1/\mu_{ijl}=\infty$ and thus \eqref{e-A3} gives
$\lambda_{ijkl}=0$ in agreement with \eqref{e-A2}. The cases where
any of $\mu_{ikl},\mu_{jkl},\mu_{kij},\mu_{lij}$ are 0 or $\infty$
are treated in the same way. (d) $\mu_{ikl}=1$. Then $\mu_{kil}=0$
and thus $\mu_{kij}\mu_{kjl}=0$. The case $\mu_{kij}=0$ is covered
by (c) so let $\mu_{kjl}=0$ whence $\mu_{jkl}=1$. Eq.~\eqref{e-A2}
gives then $\lambda_{ijkl}=1$. By $(\mu3)$,
$\mu_{ikj}=\mu_{ikl}\mu_{ilj}=1\cdot\mu_{ilj}$, $(\mu1)$ implies
$\mu_{kij}=\mu_{lij}$ and thus also \eqref{e-A3} gives
$\lambda_{ijkl}=1$. The remaining case $\mu_{kij}=\mu_{lij}=1$ is
treated in the same way by exchanging $i,j$ and $k,l$.

Thus $g$ is a well-defined morphism $Y_n\to \bar M_{0,n+1}$ and by
construction $g\circ f$ is the identity.

It remains to show that $\lambda=g(\mu)$ obeys the relations
$(\lambda1)$--$(\lambda3)$. The first relation $(\lambda1)$ is
obviously satisfied. The relation $(\lambda2)$ is satisfied by
construction if one of $i,j,k,l$ is equal to $n+1$. Also by
construction we have
\begin{equation}\label{e-A5}
\lambda_{ijkl}=\lambda_{klij},
\end{equation}
for all distinct $i,j,k,l$. If $(i,j,k,l)\in V_n$ and
$\mu_{jli},\mu_{kli}\notin\{0,\infty,1\}$, $(\mu1)$--$(\mu3)$ imply
\begin{gather*}
\lambda_{jkli}=\frac{\mu_{jli}}{\mu_{kli}}
=\frac{1-\mu_{lji}}{\mu_{kli}}
=\frac{1-\mu_{ljk}\mu_{lki}}{\mu_{kli}}
=\frac{1-(1-\mu_{jlk})\mu_{lki}}{1-\mu_{lki}}
\\
=1+\frac{\mu_{jlk}\mu_{lki}}{1-\mu_{lki}}
=1-\frac{\mu_{ikl}}{\mu_{jkl}}=1-\lambda_{ijkl}.
\end{gather*}
We now consider the degenerate cases. (a) If $\mu_{kli}=0$ then
$\mu_{lki}=1$ and thus $\mu_{lji}=\mu_{ljk}\mu_{lki}=\mu_{ljk}$ and
therefore also $\mu_{jli}=\mu_{jlk}$. Either $\mu_{jli}\neq 0$ so
that $\lambda_{jkli}=\infty$ and
$\lambda_{ijkl}={\mu_{ikl}}/{\mu_{jkl}}=\infty\cdot\mu_{jlk}=\infty\cdot\mu_{jli}=\infty$
proving the identity; or $\mu_{jli}=0$ so that the second formula
\eqref{e-A3} applies and we have
\[
\lambda_{jkli}
=\frac{\mu_{ljk}}{\mu_{ijk}}
=\frac{\mu_{lji}}{\mu_{ijk}}
=\frac1{\mu_{ijk}}=\mu_{ikj}.
\]
On the other hand, since $\mu_{lij}=1/(1-\mu_{jli})=1$,
\[
1-\lambda_{ijkl}=1-\frac{\mu_{kij}}{\mu_{lij}}=1-\mu_{kij}=\mu_{ikj},
\]
proving the claim. (b) If $\mu_{kli}=1$ then
$\mu_{ikl}=0=\mu_{lki}$. Then either $\mu_{jkl}=0$ and (a) with
permuted indices gives $\lambda_{ijkl}=1-\lambda_{lijk}$ which with
\eqref{e-A5} implies the claim; or $\mu_{jkl}\neq 0$ and
$\lambda_{ijkl}=\mu_{ikl}/\mu_{jkl}=0$. In this case
$\lambda_{jkli}=\mu_{jli}/1=1-\mu_{lji}=1-\mu_{ljk}\mu_{lki}=1-\mu_{ljk}\cdot
0=1$, since $\mu_{ljk}=1-1/\mu_{jkl}\neq\infty$. (c) If
$\mu_{kli}=\infty$ then $\mu_{ilk}=1$ and we are in case (b) up to
permutation of $i$ and $k$, so that we get
$\lambda_{jilk}=1-\lambda_{kjil}$, which reduces to the claim by
using $(\lambda1)$ and \eqref{e-A5}. Thus the cases where the
denominator $\mu_{kli}$ belongs to $\{0,\infty,1\}$ are covered by
(a)--(c). (d) The case where the numerator $\mu_{jli}$ is in
$\{0,\infty,1\}$ is reduced to the previous case by the substitution
$i\leftrightarrow k$, $j\leftrightarrow l$. Indeed (a)--(c) give
$\lambda_{lijk}=1-\lambda_{klij}$, which reduces to $(\lambda2)$ by
applying \eqref{e-A5}. This completes the proof of $(\lambda2)$.

Finally, if $\mu$ is generic, $(\lambda3)$ follows from $(\mu3)$:
\begin{equation}\label{e-A6}
\lambda_{ijkl}\lambda_{ijlm}=\frac{\mu_{ikl}}{\mu_{jkl}}\frac{\mu_{ilm}}{\mu_{jlm}}
=\lambda_{ijkm}.
\end{equation}
This formula applies more generally if one or both $\lambda_{ijkl}$
and $\lambda_{ijlm}$ are given by \eqref{e-A2}: the only tricky case
is if the left-hand side of \eqref{e-A6} is $0\cdot\infty$, but in
this case there is nothing to prove (see the footnote on page
\pageref{f-1}.) If both factors are given by \eqref{e-A3}, we have
(trivially):
\[
\lambda_{ijkl}\lambda_{ijlm}=\frac{\mu_{kij}}{\mu_{lij}}\frac{\mu_{lij}}{\mu_{mij}}
=\lambda_{ijkm}.
\]
The remaining case is when for one factor, say $\lambda_{ijkl}$,
\eqref{e-A3} is not defined and for the other, say $\lambda_{ijlm}$,
\eqref{e-A2} is not defined. Then
\[
\lambda_{ijkl}=\frac{\mu_{ikl}}{\mu_{jkl}},\qquad\lambda_{ijlm}=\frac{\mu_{lij}}{\mu_{mij}}.
\]
If \eqref{e-A3} for $\lambda_{ijkl}$ is $0/0$, i.e., $\mu_{lij}=0$
and $\mu_{kij}=0$, we have that $\mu_{mij}\neq 0$ (since
$\lambda_{ijlm}$ is assumed to be defined by \eqref{e-A3}) and
$\lambda_{ijlm}=0/\mu_{mij}=0$. Also
$\lambda_{ijkm}=\mu_{kij}/\mu_{mij}$ is defined and equal to zero
and $(\lambda3)$ is obeyed. Similarly, in the case where
\eqref{e-A3} is $\lambda_{ijkl}=\infty/\infty$, $(\lambda3)$ is
obeyed since $\lambda_{ijlm}=\lambda_{ijkm}=\infty$.
\end{proof}

\end{document}